\documentclass[reqno, oneside, 12pt]{amsart}

\usepackage[utf8]{inputenc}
\usepackage{hyperref}
\hypersetup{
    colorlinks=true,
    linkcolor=blue,
    filecolor=magenta,      
    urlcolor=cyan,
}
\usepackage{amssymb, xcolor,graphicx}
\usepackage{fullpage}
\usepackage{enumitem}
\usepackage{setspace} 
\usepackage{cuted}
\usepackage{wrapfig}

\title{Examining the Modeling Framework of Crime Hotspot Models in Predictive Policing} 

\author{Heidi Goodson}
\address{Department of Mathematics, Brooklyn College, City University of New York; 2900 Bedford Avenue, Brooklyn, NY 11210}
\email{heidi.goodson@brooklyn.cuny.edu}

\author{Alanna Hoyer-Leitzel}
\address{Department of Mathematics and Statistics, Mount Holyoke College; 50 College Street,
South Hadley, Massachusetts 01075}
\email{ahoyerle@mtholyoke.edu}

\begin{document}

\maketitle

\begin{abstract}
    Predictive policing has its roots in crime hotspot modeling. In this paper we give an example of what goes into mathematical crime hot spot modeling and show that the modeling assumptions perpetuate systemic racism in policing. The goal of this paper is to raise objections to this field of research, not on its mathematical merit, but on the scope of the problem formation.
\end{abstract}
 
\section*{Introduction}

We began writing this in the summer of 2020 after yet another Black person, George Floyd, was murdered at the hands of a member of the Minneapolis Police Department on May 25. This tragedy and the protests against police brutality and racial injustice that followed motivated many people to start or refocus their work on regular anti-racist actions. As mathematicians, one of our responses was to write this article. Our objective is to expose the terms, assumptions, and consequences of the theoretical framework that motivates predictive policing. Through this critical analysis we clarify how the model itself encodes systemic racism.

Predictive policing uses models, data, and algorithms in an attempt to predict and deter crime by optimizing police departments’ allocations of resources, including officers’ labor. At present, predictive policing research  primarily focuses on optimizing police officer deployment to crime hotspots and operates under the assumption that an officer's presence deters crime. Academic math research has led to some of the predictive policing software used by U.S. police departments \cite{PrePolPatent}. 
Many people have written about how the data and  algorithms of predictive policing software encode biases and racism (see, for example, \cite{Lum2016}, \cite{ONeil}, and \cite{PolicingOpenLetter}). Here we focus instead on the fundamental assumptions of the models justifying predictive policing.

The mathematical modeling process can be simplified into three major steps: (1) making simplifying assumptions and building a model, (2) solving the mathematical problem that was created, and (3) interpreting  and assessing the results. We are focusing on the first step in this paper, and we believe that this step is as much the responsibility of  mathematicians as the following two steps. For this reason, this article may feel more expository than what one usually sees in a review paper. 
We aim to show that the underlying agent-based and partial differential equations (PDE) models behind predictive policing are built on biases and racist policies.

\section*{Theoretical Predictive Policing Models}\label{Sec:PPModel}

In this section, we develop a detailed analysis of mathematical research justifying the development of predictive policing. After a broad search of pertinent literature, we opted to focus on twenty papers in the field of crime hotspot research, as the idea of ``crime hotspots" is fundamental to predictive policing.\footnote{We found seventeen papers using the query ``Anywhere=(crime hotspot)" on MathSciNet. We included three other papers that were highly cited by those in the search. We are happy to provide a list upon request.} We examined these papers for common assumptions and modeling techniques.  All but one of the twenty papers used agent-based and the continuum limit PDE models to create crime hotspots. At least five of them included police as part of the model. 

We chose to frame our discussion in this article around the paper ``Cops  on  the  Dots  in  a  Mathematical Model of Urban Crime and Police Response'' \cite{CopsonDots}. ``Cops on Dots" is representative of how predictive policing models are developed and presented in the mathematical literature because it extends an agent-based model to a PDE continuum model. It also clearly states that the goal of predictive policing is optimal police deployment for crime deterrence. 
While not every paper states this as clearly, it is the overarching goal of crime hotspot modeling.
It does differ slightly from other models in our search by modeling police through a deterrence to crime in the attractiveness field rather than as separate actors in the agent-based model (see points 6 - 10 of the model description).  Significantly, this paper and most of the others in our search are purely mathematical and do not employ data or algorithms (other than the algorithms needed to implement the numerical solutions for the model equations) in order to focus on the theoretical basis of predictive policing, rather than on issues with application.

\subsection*{Description of the Model}

There are two main steps to modeling crime and policing in \cite{CopsonDots}.  First, crime is modeled using an agent-based model, first published in ``A  Statistical  Model  of  Criminal  Behavior" \cite{Shortetal2008}. An agent-based model is a type of 
model based on the main concepts of interacting entities (agents/criminals)  and events (crime). 
In the second step, the discrete agent-based model in \cite{CopsonDots} is converted to a continuous model and police presence is added as a variable   optimized to give the minimum crime density in the space. The development of the model in \cite{CopsonDots} follows this outline:
    \begin{enumerate}\setlength\itemsep{0em}
        \item   A lattice is defined for the space $\Omega$ on which the agents in the agent-based model may move around. In each time step, an agent has some probability of moving to a different lattice point. 
        \item In the agent-based model, every agent is called a {\it criminal}. At every time step, at the locations of each criminal\footnote{In the body of \cite{CopsonDots}, non-police agents are called ``criminals'' 75\% of the time.  We find the use of the word ``criminal" problematic, since an agent may not actually commit a crime, they only have some nonzero probability of doing so. The presumption of guilt here is part of how such models instill systemic, to-some-extent-unconscious bias. Unfortunately, we will use that language here because that is used in the paper.}, there is some nonzero probability that an event occurs based on an attractiveness field $A(x,t)$. These events are called {\it crimes}.  The attractiveness field assigns values to the lattice nodes;  a higher value means that it is more likely that a criminal will commit a crime at that node. 
        \item When a criminal commits a crime, they disappear from the lattice. If a criminal is removed from the system, there is some probability that more criminals will enter the system.
        \item The attractiveness field is updated dynamically. If a crime is committed at a node, the attractiveness of that node and nearby nodes increases. This creates {\it crime  hotspots}, locations in $\Omega$ where more crimes  occur than in nearby locations.
        \item In the limit, as the lattice point spacing and time step length go to zero, the agent-based model can be written as a system of continuous PDEs with two equations: one PDE for the attractiveness field, $\partial A/\partial t$, and another PDE for the expected value of the density of criminals in the space, $\partial \rho /\partial t$.  Stationary solutions of this PDE system generate hotspots like the agent-based model. 
        \item  {\it Cops} are added to the PDE system with a spatial and time dependent function $\kappa(x,t)$, the number of cops deployed to point $x$ at time $t$. This function $\kappa$ is composed with a deterrence function $d:[0,\infty)\to(0,1]$, which inputs some number of cops and outputs a scaling factor which is multiplied against the product $\rho A$, 
        a term that describes the amount of crime and which appears in both $\partial A/\partial t$
        and $\partial \rho /\partial t$. 
        \item The deterrence function satisfies, among other conditions:  (a)         When $d(\kappa(x,t))=1$ for $x\in\Omega$ and time $t$, the system is the same as if there were no cops,         and (b)  $\lim_{\kappa\to\infty}d(\kappa)=0$, i.e., the police ``can achieve a target deterrence level if they deploy enough resources''  \cite{CopsonDots}.
        \item Cop deployment is modeled by solving an optimization problem.  The goal of the optimization is to minimize crime in the system with the constraint that the cops have a finite amount of resources, i.e. minimize $\int_{\Omega}d(\kappa(x,t))\rho(x,t)A(x,t)\,dx$ with the constraint that $\int_{\Omega}\kappa(x,t)\,dx=K$ for some fixed $K\geq0$ (among other technical constraints).
      \item   The effect of police deterrence varies with the choice of $K$ (the total number of cops). Smaller $K$ reduces the intensity of hotspots but makes the spots bigger. Increasing $K$ more will cause the hotpots to merge and form ``worm-like" features in $\Omega$. High values of $K$ will cause total hotspot suppression so that $\rho$ and $A$ are uniform on the domain.
        \item  Police deterrence in the attractiveness field is then incorporated back into the agent-based model with similar results to the PDE model.
        \end{enumerate}

\section*{Racism in the Modeling Assumptions}\label{ModelProblems}

People often think that racism is about intentional interpersonal interactions that incite violence, physically or emotionally harming a person on the basis of race.  But scholars and activists working in areas of racial justice often use racism to describe a more systemic issue. Here we use {\it racism} to mean systemic racism, where policies or institutions create disparate outcomes for social groups, and race is one clear axis on which disparate outcomes can be measured. Systemic racism is reinforced by systems of power; it ``is a machine that runs whether we pull the levers or not, and by just letting it be, we are responsible for what it produces" \cite{Oluo}.  While we focus on anti-Black racism, systemic racism also affects Latinx people, Native Americans, immigrants, and other marginalized and minority groups.

In modeling, it is common  to assume a homogeneous population across gender, race, age, or other factors. This creation of an abstract population devoid of demographics is akin to white privilege. White privilege refers to when white people may lack -- or even be unaware of -- experiences common to marginalized people of color. ``Abstraction privilege" in math plays a large role in how we create and evaluate models. In what follows, we discuss four ways that the modeling that motivates crime hotspot and police deployment optimization problems ignores or obscures race. By ignoring race, the model ignores the effects of racism and, thus, perpetuates it.\\

\noindent {\bf Crime and criminality are racialized concepts in the United States, but this is not acknowledged in the model context.} The United States has a history of constructing Black people (as well as Indigenous people and many other people of color) as criminals even in the absence of suspected or actual crimes.  Black people were labeled ``criminals" in racist justifications of slavery prior to the Civil War.  ``From their arrival around 1619, African people had illegally resisted slavery. They had thus been stamped from the beginning as criminals" \cite{Kendi}.  Policing in the United States developed out of slave patrols in the U.S. South as well as Night Watches in the North, both of which functioned to control impoverished groups, particularly Black people and Native Americans  \cite{Vitale}. 
Thus, ideas of ``crime” in the U.S. are historically tied to white colonists’ and landowners' desires to maintain power over other populations for the purpose of extracting value from their labor.   

The practice of labeling Black people ``criminals” permeates contemporary society. Sales clerks surveil Black customers, white people avoid Black people on the street and view them as more likely to commit crimes. Media portrayals of race and crime exacerbate these biases \cite{Oliver2003}. This leads directly to interactions with police and the criminal justice system.  For example, a 2009 study of stop-and-frisk data in New York City by the Center for Constitutional Rights found that the rate at which Black and Latino New Yorkers were stopped and frisked by NYPD officers was significantly disproportionately higher than for white New Yorkers, and that these rates do not correspond to rates of arrests or summons \cite{CCR2009}. This criminalization of Black people starts early as, tragically, Black children in the U.S. are disproportionately funneled into the prison system via the ``school-to-prison pipeline”. In 2011 - 2012 Black students represented 27 percent of students referred to law enforcement and 31 percent of students subjected to a school-related arrest, though they only made up 16 percent of student enrollment in the U.S.  \cite{Vitale}.  Once in the prison system, Black people are subjected to harsher sentences and parole decisions: a 2016 article found that Northpointe Inc.’s COMPAS recidivism algorithm falsely labeled Black defendants as future criminals at almost twice the rate as white defendants \cite{Propublica},  resulting in longer sentences and lower likelihood of parole.

Crime hotspot and police deployment models do not acknowledge the racialized nature of notions of crime and criminality in the United States. While it may be convenient to create a model with a homogeneous population that does not explicitly take race into account, there should be acknowledgment that this assumption affects potential implementations of policing based on the results. The absence of such acknowledgement perpetuates a facade of non-bias for models that are implicitly biased.\\

\noindent\textbf{The way the attractiveness field updates is justified by a controversial policing theory.}
The definition of the attractiveness field in many crime hotspot models (see point 4 in the description of the model) is from a paper by Short et al. \cite{Shortetal2008}, where the authors cite Broken Windows Theory.  Broken Windows Theory comes from an Atlantic Monthly article \cite{Atlantic1982} suggesting that  a neighborhood that tolerates disorder in the form of ``small crimes" (loitering, pan handling, homelessness, jumping subway turnstiles) will also tolerate more serious crimes (burglary and violence). This kind of crime concentration is often designated a \textit{crime hotspot}. In the models we examined, as `disorder’ in the form of crimes accumulates at a given lattice point, the model predicts that further crimes will be drawn to that point.

  \begin{figure}
  \begin{center}
  \includegraphics[width=.9\textwidth]{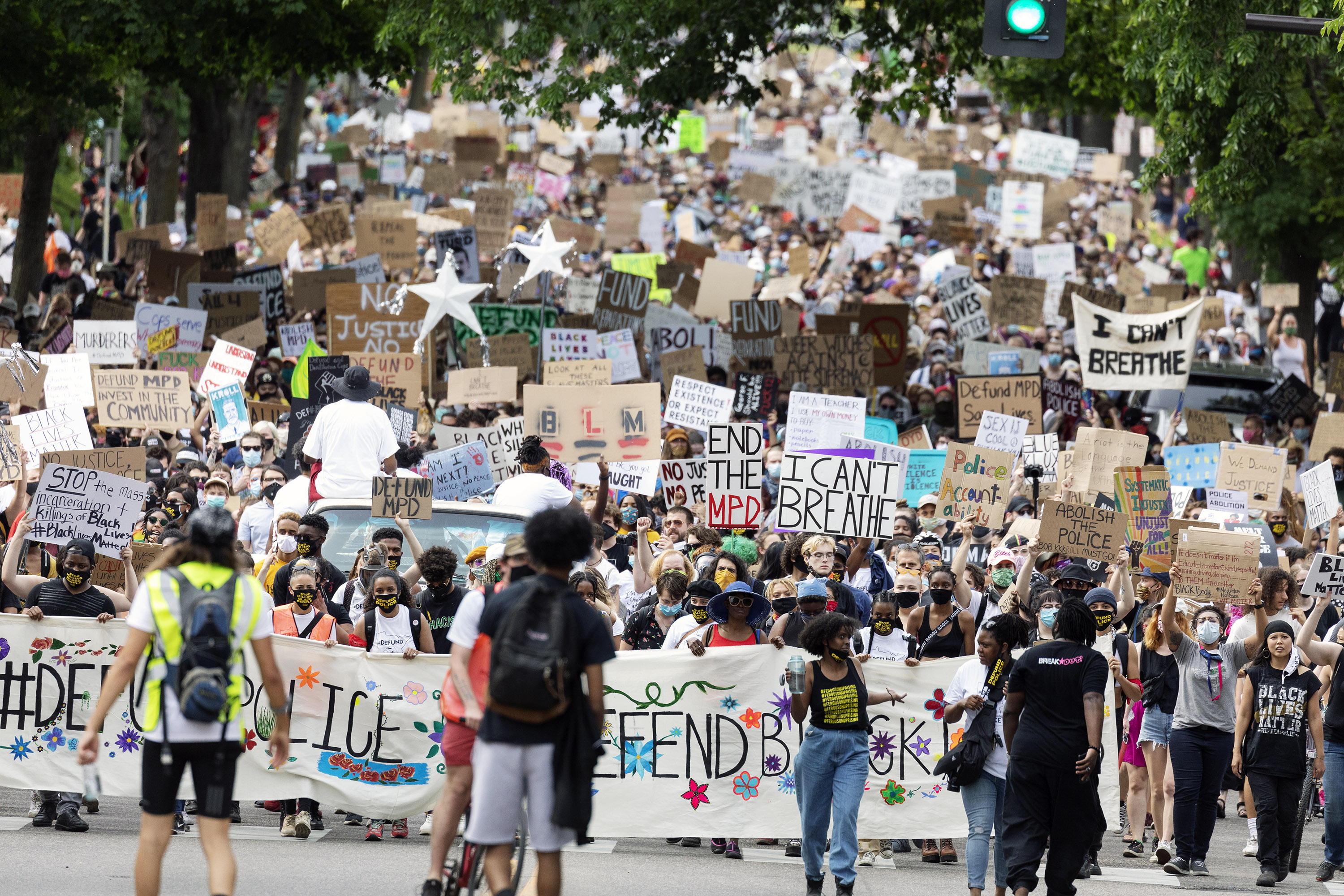}\\
  \caption{Defund the MPD Rally in the wake of George Floyd's murder. Minneapolis, MN. June 6, 2020. Photo credit: Ben Hovland.}
  \end{center}
  \end{figure}
In reality, our biases play a role in what we view as disorder and in what suspected crimes law enforcement officers investigate. A 2004 article by Sampson and Raudensbush concludes: ``as the concentration of minority groups and poverty increases, residents of all races perceive heightened disorder"  \cite{Sampson}. Historically, behaviors that do not match those of people who are wealthy and/or white  are deemed disorderly or criminal. There has been  much discussion (see, for example, \cite{GreeneZeroTolerance, Vitale}) about how Broken Windows Theory, as applied in the United States, has led to practices, such as zero-tolerance and stop-and-frisk policing, which have disproportionately targeted people of color and unhoused people. Our main objection is that Broken Windows Theory fails to address  the complicated relationships among perceived disorder, poverty and Blackness in the United States. 

The authors of crime hotspot modeling papers do not acknowledge the controversies of Broken Windows Theory. Hence, the model’s role in perpetuating racialization of policing is further obscured. \\

\noindent \textbf{While the model treats location as neutral, geography is racially segregated.} 
In the models  developed in \cite{CopsonDots} and \cite{Shortetal2008}, the attractiveness field is treated as neutral in the sense that underlying attributes about the location of the agent are not taken into account. However, in many ways, indicating a person's location is (at least probabilistically) indicating the race of that person. In hundreds of cities across the United States, residential geography is heavily shaped by redlining, an explicitly racist housing policy from the early 20th century that used race to guide real estate practices (see, for example, \cite{HOLCreport}, \cite[Chapter 2]{AmericanApartheid}, and \cite[Chapter 28]{Kendi}). Through redlining,  white home-ownership and suburban development were subsidized while resources were withheld from people of color. This created the racialized patterns of urban segregation, disinvestment, and poverty that still exist today. While redlining was made illegal by the Fair Housing Act in 1968, neighborhoods that were redlined in the 1930s have higher rates of poverty even today. According to the National Community Reinvestment Coalition, 3 out of every 4 neighborhoods in the United States that were redlined in the 1930s are  low-to-moderate income today, and 2 out of every 3 are predominantly populated by people of color  \cite{HOLCreport}.  

Predictive policing models are developed to potentially influence public policy. However when building the assumptions of the model the authors did not take into account issues like redlining, which contribute to the ways that there are racial geographic patterns. \\

\noindent{\bf The societal impacts of police behavior are not reconciled with the simplicity of the model.}  
The model assumes that the only effect of police presence is to deter crime (see point 6 of the model description).
However, when cops are actually deployed, they  interact with people and have the power to do harm. Indeed, we are writing this article specifically because George Floyd, Breonna Taylor, Tony McDade, Philando Castile, Eric Garner, Michael Brown, Tamir Rice, and so many other Black people have been murdered by police, often following interactions that should have been routine and non-criminalized but were unnecessarily escalated by the police. 

The dangers of policing are   particularly pronounced for people of color. Using national data collected between 2013 and 2018, \cite{Edwards2019} found ``that people of color face a higher likelihood of being killed by police than do white men and women, that risk peaks in young adulthood, and that men of color face a nontrivial lifetime risk of being killed by police."
In particular, they found that Black  men  are  about  2.5  times more  likely  than white  men to  be  killed  by  police.

Police interaction and surveillance affect not only people who are stopped by police but also those who live in highly-policed communities. In a 2016 article, Sewell et al. examined the mental health effects for people  living in ``aggressively policed" New York City neighborhoods \cite{Sewell2016}. Their analyses show that ``living in aggressively policed communities is of detriment to the health of male residents" and  these men are more likely to experience severe psychological distress and other mental health issues.\footnote{Most of the conclusions in \cite{Sewell2016} were about the effects on men. The authors of the study noted that over 85 percent of pedestrians stopped by NYPD were men.}

Harm done by cops is not included in the development or analysis of predictive policing models. The crime hotspot and police deployment models we examined consider only how police presence deters crime -- one small aspect of policing and police behavior.

\section*{Concluding Remarks}

The mathematical content of crime hotspot and police deployment models is developed and presented with the usual mathematical rigor and novelty one expects of published work. Often a goal of modeling is to examine what mechanisms create key processes in a system. It is common practice to build a model using simplifying assumptions in order to make analysis tractable. But one must be responsible about the influence of these models, particularly when they have influence on  public policy and decision making. We believe that crime hotspot and police deployment optimization models examine the key processes of a simplified system built to deter crime, but with assumptions that encode systemic racism through abstraction and false neutrality. 

Mathematical research developing predictive policing models as it currently exists in the literature fails to adequately consider the social and racial contexts of its application, as well as how policing perpetuates and may exacerbate systemic racism. On this basis, we believe that current models should have no role in answering the question of how police might be efficiently and ethically deployed. Additionally, using these theoretical models to examine the mechanisms of implicit bias and racism in policing, or to quantify racism in policing is misguided. Instead, we should listen to scholars, activists, and Black people when they say that there are issues with policing in the United States, particularly around race.  For these reasons, it is difficult for us to envision how to do anti-racist predictive policing research because the problems discussed above are so fundamental.

We wonder whether police deployment is the right optimization problem to consider. Policing addresses the effects of crime rather than the causes of crime: poverty and lack of social resources \cite{Vitale}. If the goal is to optimize the use of limited government and community resources, the optimization problem should focus instead on how to distribute social supports such as housing, food, and mental and physical healthcare.

The role of mathematicians in policing and the development of predictive policing is an ongoing and important conversation for us to continue (see, for example, \cite{JMCblog} and the Letters to the Editor in the September and October 2020 Notices of the American Mathematical Society). This is part of the larger national conversation about policing happening in the United States. We hope that this article contributes to this discussion by offering an exposition of predictive policing’s theoretical framework. We feel strongly that researchers and the mathematical institutions that support them need to reevaluate the criteria by which they examine the greater effects -- the true broader impacts -- of this kind of research to ensure that mathematics fights, rather than maintains, systemic racism.



\bibliographystyle{plain}
\bibliography{policing.bib}

\end{document}